\documentclass[11pt]{article} 
\usepackage{latexsym}
\usepackage{amsfonts}
\usepackage{amssymb}
\usepackage{graphicx}
\begin{document}

\title{Matrix Theory over the Complex Quaternion Algebra}
\author{Yongge   Tian \\
Department of Mathematics and Statistics \\
Queen's University \\
Kingston,  Ontario,  Canada K7L 3N6 \\
{\tt email:ytian@mast.queensu.ca}}
\date{}
\maketitle 

\noindent {\bf Abstract.} \, {\small We present in this paper some fundamental tools for developing matrix analysis over the complex quaternion algebra. As applications, we consider generalized inverses, eigenvalues and eigenvectors, similarity, determinants of complex quaternion matrices, 
and so on. 

\bigskip

 }  \\

\bigskip

\noindent {\Large  {\bf 1. \ Introduction}} \\

\noindent The complex quaternion algebra (biquaternion algebra) $\Bbb Q$ is 
well known as a four dimensional
vector space over the complex number field $ \Bbb C $ with its  basis
$ 1, \, e_1, \, e_2, \, e_3 $ satisfying the  multiplication laws
$$ 
e^2_1 =  e^2_2 = e^2_3 = -1, \ \ \   e_1 e_2 e_3 = -1.  \eqno (1.1)  
$$ 
$$ 
 e_1 e_2 = -e_2e_1 = e_3,   \ \ e_2 e_3 = -e_3e_2 = e_1, \ \
 e_3 e_1 = -e_3e_3 = e_2,
 \eqno (1.2)  
$$  
and 1 acting as unity element. In that case, any element in $ \Bbb Q $ can
be written  as
$$ 
a = a_0 + a_1e_1 + a_2e_2 +a_3 e_3, \eqno (1.3) 
$$ 
where $a_0$---$a_3 \in \Bbb C$. According to this definition, real numbers,
complex numbers, and real quaternions all can be regarded as the special
cases of complex quaternions. A well-known fundamental fact on the complex
quaternion algebra $\Bbb Q$ (see, e. g., \cite{4, 6, 7}) is
that it is algebraically isomorphic to the $2 \times 2$ total matrix
algebra  $ \Bbb C^{2 \times 2}$ through the bijective map $ \psi :
\Bbb Q \longrightarrow \Bbb C^{2 \times 2}$ satisfying
$$ 
 \psi(1) = \left[ \begin{array}{cc} 1 & 0 \\ 0 & 1 \end{array} \right], \ \
 \ \psi(e_1) = \left[ \begin{array}{cc} i & 0 \\ 0 & -i \end{array} \right],
 \ \ \ \psi(e_2) = \left[ \begin{array}{cc} 0 & -1 \\ 1 & 0 \end{array}
 \right], \ \ \ \psi(e_3) = \left[ \begin{array}{cc} 0 &  -i \\ -i & 0
 \end{array}  \right].
$$
These four matrices are well-known as Pauli matrices. Based on this map,
every element
$a = a_0 + a_1e_1 + a_2e_2 +a_3 e_3 \in \Bbb Q$ has a faithful  complex
matrix representation  as follows
$$
 \psi(a):= \left[ \begin{array}{cc}  a_0 + a_1 i & -( \, a_2 + a_3i \, ) \\
 a_2 - a_3 i & a_0 - a_1 i  \end{array}  \right] \in \Bbb C^{ 2 \times 2},
\eqno (1.4)
$$  
In this article, we shall reveal a deeper relationship between $ a $ and
$ \psi(a)$, which can simply be stated that there is an independent
invertible matrix $Q$ of size 2 over $ \Bbb Q$ such that all $ a \in
\Bbb Q$ satisfy the following universal similarity factorization equality
$$
Q^{-1}{\rm diag}( \, a, \ a \, )Q =  \psi (a), 
$$         
where $ Q $ has no relation with the expression of $ a $. Moreover we also
extend this equality to all $ m \times n $ matrices over $ \Bbb Q$. On the
basis of these results, we shall consider several basic problems related to
complex quaternion matrices, such as, generalized inverse, eigenvalues and
eigenvectors, similarity, and determinant of complex quaternion matrices. 

\medskip

Some known terminology on complex quaternions are listed below (see, e.g.,
\cite{6}). For $ a = a_0 + a_1e_1  + a_2e_2  + a_3 e_3 \in \Bbb Q,$
 the  {\em dual quaternion } of $ a $ is
$$  
\overline{a} = a_0 - a_1e_1  - a_2e_2  - a_3 e_3; \eqno (1.5) 
$$ 
the {\em complex conjugate} of $ a $ is 
$$ 
a^* = \overline{a_0} +  \overline{a_1}e_1  +  \overline{a_2}e_2  +
\overline{ a_3} e_3; \eqno (1.6)
$$ 
the {\em Hermitian conjugate} of $ a $ is 
$$ 
a^{\dagger} = ( \overline{a} )^*  =  \overline{a_0} -
\overline{a_1}e_1  -  \overline{a_2}e_2  - \overline{ a_3} e_3; \eqno (1.7)
$$ 
the {\em weak norm} of $ a $ is   
$$ 
n(a) = a_0^2 + a_1^2+ a_2^2 +a_3^2 . \eqno (1.8)
$$ 
A quaternion $ a \in \Bbb Q $ is said to be {\em real} if $a^*  = a$, 
 to be {\em pure imaginary} if $a^*  = -a$, to be  {\em scalar} if $ \overline{a} = a$, to be   {\em Hermitian} if $a^{\dagger} = a$.
 
\medskip

For any  $ A = (a_{st}) \in  \Bbb Q^{m \times n}$, the {\em dual}  of
$ A $ is $\overline{ A } = ( \overline{ a_{ts} } ) \in
 \Bbb Q^{n \times m}$; the {\em Hermitian conjugate} of $ A $ is
 $ A^{\dagger} = (a_{ts}^{\dagger}) \in \Bbb Q^{n \times m }$.
 A square matrix $ A $ is said to be {\em  self-dual} if
 $\overline{ A } = A$, it is  {\em Hermitian} if $ A^{\dagger} = A$,  it is
  {\em unitary} if $AA^{\dagger}= A^{\dagger} A = I$, the indentity matrix,
   it is {\em invertible} if there is a matrix $ B $ over ${\cal  Q}$ such
   that $AB = BA = I$.  

\medskip

Some known basic properties on complex quaternions and matrices of complex
quaternions are listed below. 

\medskip

\noindent 
{\bf Lemma 1.1}\cite{1}\cite{6}. \, {\em Let $ a, \, b \in \Bbb Q$ be
given. Then

{\rm (a)} \ $ \overline {\overline{a}} = a,  \ \ \ (a^*)^* = a,  \ \ \
(a^{\dagger})^{\dagger}= a. $
  
{\rm (b)}  \ $ \overline{ a + b }=  \overline{a} +  \overline{b},
\ \ \ \ \ ( a + b )^* = a^* + b^*,  \ \ \ \ \ ( a + b )^{\dagger}
= a^{\dagger} + b^{\dagger}.$

{\rm (c)} \ $ \overline{ab}=\overline{b}\overline{a}, \ \ \ \ \ (ab)^* 
= a^*b^*, \ \  \ \ \ (ab)^{\dagger} = b^{\dagger}a^{\dagger}.$  

 {\rm (d)} \ $ a\overline{a} = \overline{a}a = n(a) = n(\overline{a});  $   

{\rm (e)} \ $ a $ is invertible if and only if $ \ n( a ) \neq 0,$
in that case $ a^{-1} =  n^{-1}( a ) \overline{a}. $ }

\medskip

\noindent 
{\bf Lemma 1.2}\cite{6}. \, { \em Let $ A \in \Bbb Q^{m \times n}, \, B
\in \Bbb Q^{n \times p }$ be given. Then

{\rm (a)} \ $ \overline {\overline{A}} = A,   \ \ \
(A^{\dagger})^{\dagger}= A.$

 {\rm (b)} \  $\overline{AB} = \overline{B} \overline{A} ,
 \  \ \ \ \ \ \ ( AB )^{\dagger} = B^{\dagger}A^{\dagger}.$

{\rm (c)} \ $( AB )^{-1} = B^{-1}A^{-1},$ \ if $ A $ and $ B $ are
invertible.

{\rm (d)}  \ $ (\overline{ A})^{-1} = \overline{( A^{-1})},
\ \ \ \ \  ( A^{\dagger})^{-1} =  ( A^{-1})^{\dagger},$ if $ A$ is
invertible.} \\

\noindent{\Large {\bf 2. A universal similarity factorization equality over
complex quaternion algebra}} \\
  
We first present a general result on the universal similarity factorization
of elements over $ 2 \times 2 $  total matrix algebra. 

\medskip

\noindent 
{\bf Lemma 2.1.} \, {\em  Let $M_2(\Bbb F )$ be the  $2 \times 2 $ total
matrix algebra  over an arbitrary field $\Bbb F$ with its basis
$ e_{11}, \, e_{12}, \,  e_{21}$ and $e_{22} $ satisfying the following
multiplication rules
$$
e_{st} e_{pq} = \left\{ \begin{array}{c} e_{sq},  \ \ \ t = p \\ 0,  \ \ \ \ \ t \neq p \end{array} 
\right. , \ \ \ \ s, \, t, \, p, \, q = 1, \, 2. \eqno (2.1) 
$$  
Then for any  $a = a_{11}e_{11}  + a_{12}e_{12}  +  a_{21}e_{21}  +
a_{22}e_{22} \in M_2({\cal F}),$ where $a_{st} \in {\cal F}$, the
corresponding diagonal matrix $ {\rm diag}( \, a, \ a \, ) $ satisfies
the following  universal similarity factorization equality
$$
Q \left[ \begin{array}{rr}  a  & 0   \\ 0 &  a \end{array} 
\right] Q^{-1}= \left[ \begin{array}{cc}  a_{11} & a_{12}  \\a_{21} & a_{22}
\end{array} \right] \in {\cal F }^{2 \times 2},
  \eqno (2.2) 
$$  
where $Q$ has the independent form 
$$ 
Q= Q^{-1} = \left[ \begin{array}{cc} e_{11} & e_{21}  \\  e_{12} &  e_{22}
\end{array} \right].  \eqno (2.3)
$$ }
\noindent {\bf Proof.} \,  According to Eq.(2.1), it is easy to verify that
the unity element in $M_2({\cal F})$ is $ e =  e_{11} + e_{22}$. In that case,
 the matrix $ Q $ in Eq.(2.3) satisfies
\begin{eqnarray*} 
Q^2 = \left[ \begin{array}{cc} e_{11} & e_{21} \\ e_{12} & e_{22}
\end{array} \right]  =
\left[ \begin{array}{cc} e_{11}^2 + e_{21}e_{12} & e_{11} e_{21}+ e_{21}
e_{22}  \\  e_{12} e_{11} + e_{22}e_{12} &
 e_{12} e_{21} + e_{22}^2  \end{array} \right] =  \left[ \begin{array}{cc}
  e_{11} + e_{22} & 0  \\  0 & e_{11} +  e_{22} \end{array} \right] = e I_2 ,
\end{eqnarray*} 
which implies that $ Q $ is invertible over $M_2({\cal F})$ and $Q= Q^{-1}$.
 Next multiplying the three matrices in the left-hand side of Eq.(2.2) yields
  the right-hand side of Eq.(2.2).   \qquad $ \Box  $ 

\medskip

\noindent  
{\bf Theorem 2.2.} \, { \em Let $ a = a_0 + a_1e_1 + a_2e_2 +a_3 e_3 \in
 \Bbb Q$ be given. Then the diagonal matrix $ {\rm diag }( \, a, \ a \,
  )$ satisfies the following universal  factorization  similarity  equality
$$
Q \left[ \begin{array}{rr}  a  & 0   \\ 0 &  a \end{array} 
\right] Q= \left[ \begin{array}{cc}  a_0 + a_1 i & -( \, a_2 + a_3i \, ) \\
 a_2 - a_3 i &
 a_0 - a_1 i  \end{array}
 \right] = \psi(a) \in \Bbb C^{ 2 \times 2},   \eqno (2.4) 
$$ 
where $Q$ is an unitary matrix over $\Bbb Q$ 
$$ 
Q= A^{-1} = Q^{\dagger} = \frac{1}{2} \left[ \begin{array}{cc} 1 - ie_1 &
 e_2 + ie_3 \\ -e_2 + i e_3 & 1 + ie_1  \end{array} \right]. \eqno (2.5)
$$ }
{\bf Proof.} \, According to Eq.(2.4), we choose a new basis for $\Bbb Q$ as
 follows
$$ 
e_{11} = \frac{1}{2} \left( \, 1 - ie_1 \, \right),  \ \ e_{12} = \frac{1}{2}
 \left( \, - e_2 + ie_3  \, \right), \ \ e_{21} = \frac{1}{2} \left( \, e_2
  + ie_3 \,\right ),  \ \  e_{22} = \frac{1}{2} \left( \, 1 + ie_2
   \,\right).  \eqno (2.6)
$$ 
Then it is not difficult to verify that the above basis satisfies the
 multiplication rules in Eq.(2.1). Under this new basis, any element
  $ a = a_0 + a_1e_1 + a_2e_2 +a_3 e_3 \in \Bbb Q$ can be expressed as
$$ 
a =  ( \,a_0 + ia_1 \, )e_{11}+ ( \,-a_2 - ia_3 \, )e_{12} +
  ( \,a_2 - ia_3 \,)e_{21} + ( \,a_0 - ia_1 \,)e_{22}. \eqno (2.7)
$$
Substituting Eqs.(2.6) and (2.7) into Eq.(2.2), we obtain Eqs.(2.4) and
(2.5). \qquad $ \Box $ 

\medskip
  
The equality in Eq.(2.4) can also equivalently be expressed as 
$$
 Q \left[ \begin{array}{rr}  a_{11}  & a_{12}  \\ a_{21} &  a_{22}
  \end{array}  \right] Q= \left[ \begin{array}{rr}  a  & 0   \\ 0 &  a
  \end{array}  \right] \in \Bbb Q^{ 2 \times 2}, \eqno(2.8)
$$ 
where $ a_{st} \in  \Bbb C $ is arbitrary, $ Q$ is as in Eq.(2.5), and
 $a$ has the  form
$$ 
a = \frac{1}{2}( \, a_{11} +  a_{22} \, ) +  \frac{1}{2}( \, a_{22} -
a_{11} \, )ie_1 +  \frac{1}{2}( \, a_{21} -  a_{12} \, )e_2 +
 \frac{1}{2}( \, a_{12} +  a_{21} \,)ie_3.  \eqno(2.9)
$$ 
This equality  shows that every $ 2 \times 2 $ complex  matrix is uniformly
similar to an diagonal matrix  with the form $ aI_2$ over the complex
quaternion algebra $ \Bbb Q$. 

\medskip

The complex quaternions and their complex matrix representations  satisfy
  the following operation properties.

\medskip
  
\noindent {\bf Theorem 2.3}\, { \em Let $ a = a_0 + a_1e_1 + a_2e_2 +a_3 e_3,
 \,  b \in \Bbb Q, \,  \lambda \in \Bbb C$ be given.  Then
 
{\rm (a)} \ $ a=b \Longleftrightarrow  \psi(a) = \psi(b).$ 

{\rm (b)} \ $ \psi(a + b )=  \psi(a) + \psi(b),  \ \ \psi(ab)= \psi(a)
\psi(b),\ \ \psi( \lambda a)= \psi(a \lambda ) = \lambda \psi(a),  \ \
 \psi(1) = I_2.$

{\rm (c)} \ $ \psi(\overline {a}) =  \left[ \begin{array}{rr}  0  & 1
  \\ -1 &  0 \end{array} \right] \psi^T(a) \left[ \begin{array}{rr}  0
   & -1   \\ 1 &  0 \end{array}
\right].$  

{\rm (d)} \ $ \psi(a^*) =  \left[ \begin{array}{rr}  0  & 1   \\ -1 &  0
 \end{array}
\right] \overline{\psi(a)} \left[ \begin{array}{rr}  0  & -1   \\ 1 &  0
\end{array}
\right].$ 
 
{\rm (e)} \  $ \psi(a^{\dagger}) = \overline{\psi(a)}^T = \psi^*(a),$
 the conjugate transpose of $\psi(a).$
 
{\rm (f)} \ $ \det \psi(a) =  n(a) =  a_0^2 + a_1^2+ a_2^2 +a_3^2;$ 

{\rm (g)} \  $ a = \frac{1}{4}E_2 \, \psi(a) \, E_2^{\dagger},$ where
$ E_2= [ \, 1- ie_1, \ e_2 + ie_3 \, ].$

{\rm (h)}  \  $ a$ is invertible if and only if $ \psi(a) $ is invertible,
 in that case, $\psi(a^{-1}) = \psi^{-1}(a)$ and $ a^{-1} =  \frac{1}{4}E_2
  \, \psi^{-1}(a) \, E_2^{\dagger}. $ } 

\medskip

For a noninvertible element over $ \Bbb Q $, we can define its
 Moore-Penrose inverse
as follows. 

\medskip

\noindent {\bf Definition.} \,  Let $ a \in \Bbb Q $ be given. If the
following four equations
$$
 axa = a,  \ \ \ \ \  xax = x, \ \ \ \ \ ( ax )^{\dagger}= ax, \ \ \ \ \
  ( xa )^{\dagger} = xa  \eqno (2.10)
$$ 
have a common solution $ x $, then this solution is called the Moore-Penrose
inverse of $ a $, and denoted by $ x = a^{+}$. 

\medskip

The existence and the uniqueness of the Moore-Penrose inverse of a complex
 quaternion $ a $  can be determined by its
matrix representation $ \psi(a) $ over $ \Bbb C$. In fact, according to
Theorem 2.3(a), (b) and (e), the four equations in Eq.(2.10) are equivalent
to the  following four equations over $  \Bbb C $
$$  
 \psi(a)\psi(x)\psi(a) = \psi(a), \ \ \ \ \ \  \psi(x)\psi(a)\psi(x)
  = \psi(x),
$$
$$ 
 [ \, \psi(a)\psi(x) \, ]^* = \psi(a) \psi(x) , \ \ \ \ \ \
  [ \,  \psi(x)\psi(a) \, ]^* = \psi(x) \psi(a).
$$ 
According to the complex matrix theory, the following four equations 
$$  
 \psi(a)Y \psi(a) = \psi(a), \ \ \ \ \   Y \psi(a)Y = Y,
 \ \ \ \ \ [ \, \psi( a)Y \,]^* = \psi( a)Y , \ \ \ \ \
 [\, Y \psi(a) \,]^* = Y  \psi(a).
$$ 
has a unique solution $  Y  = \psi^{+}(a)$, the Moore-Penrose inverse of
  $ \psi(a)$. Then by Eq.(2.8), it follows that there must  be a unique $ x $
  over $ \Bbb Q $ such that $ \psi( x ) = Y =  \psi^{+}(a)$,  in which
   case, this $ x $ can be expressed  as
$$
 x = \frac{1}{4}E_2  Y  E_2^{\dagger} = \frac{1}{4} E_2 \, \psi^{+}(a) \,
  E_2^{\dagger}.
$$ 
Correspondingly  this  $ x $ is the unique solution to  Eq.(2.10).  In
summary, we have the following.

\medskip

\noindent {\bf Theorem 2.4.}\, { \em Let $ a \in \Bbb Q $ be given.  Then its Moore-Penrose inverse $ a^{+} $ exists uniquely, and  satisfies the following four equalities  
$$ 
 \psi(a^{+}) = \psi^{+}(a), \ \ \ \ \ \ a^{+} = \frac{1}{4} E_2 \,
  \psi^{+}(a) \, E_2^{\dagger},
$$ 
where  $ E_2= [ \, 1 - ie_1, \ e_2 +  ie_3 \, ].$  } 

\medskip

One of the basic problems related  to complex quaternions is concerned
with similarity of two complex quaternions. As usual, two  complex
quaternions  $ a $ and $ b $ are said to be similar if there is an invertible
complex quaternion $ x $ such that $ x^{-1}ax = b $, and this is written as
 $ a \sim b $.  It is easy to  verify that the similarity mentioned here is
  an equivalence relation on  complex quaternions. 

\medskip

A simple result follows immediately from the above definition, Eqs.(2.4) and
(2.8). 

\medskip

\noindent {\bf Theorem 2.5.} \, { \em Let $ a, \,  b \in \Bbb Q $ be
 given.  Then
$$ 
a \sim b  \  \Longleftrightarrow  \ \psi(a) \sim \psi(b). \eqno (2.11) 
$$ }  

From Eq.(2.11), we easily find the following. 

\medskip

\noindent{\bf Theorem 2.6.} \, { \em Let $ a = a_0 + a_1e_1 + a_2e_2 +
a_3 e_3 \in \Bbb Q$ given with  $ a \notin  \Bbb C.$

{\rm (a)} \ If $  a_1^2+ a_2^2 +a_3^2 \neq 0, $ then $a  \ \sim
 \ a_0 + \tau(a)e_1,$ where $ \tau(a) $ is complex number satisfying
  $  \tau^2(a) =  a_1^2 + a_2^2 +a_3^2.$
  
{\rm (b)} \ If $  a_1^2+ a_2^2 +a_3^2  = 0,$ then $a  \sim
a_0 - \frac{1}{2}e_2 + \frac{1}{2}ie_3.$ }

\medskip

\noindent {\bf Proof.} \,  For any  $ a \in \Bbb Q$, the characteristic
 polynomial of its complex  matrix representation $ \psi(a) $ is
$$
 | \, \lambda I_2 -  \psi(a) \, | =  \left| \begin{array}{cc}
 \lambda -(  \, a_0 + a_1i \, ) &  a_2 + a_3i
   \\ -a_2 + a_3i  &  \lambda -( \, a_0 - a_1i \,)   \end{array} \right| =
    (  \,\lambda - a_0  \,)^2 +  a_1^2+ a_2^2 +a_3^2.
$$  
From it we immediately know that if  $  a_1^2+ a_2^2 +a_3^2 \neq 0$, then 
$$  
 \psi(a)  \ \sim  \ \left[ \begin{array}{cc}   a_0 + \tau(a)i  & 0
  \\ 0 &   a_0 - \tau(a)i  \end{array}
\right] =\psi[ \, a_0 + \tau(a)e_1 \, ],  \eqno (2.12) 
$$  
and 
if  $  a_1^2+ a_2^2 +a_3^2 =  0$, then
$$  
 \psi(a)  \ \sim  \ \left[ \begin{array}{cc}   a_0  &  1  \\ 0 &   a_0
 \end{array}
\right] =\psi\left( \, a_0  - \frac{1}{2}e_2 + \frac{1}{2}ie_3 \,\right).
 \eqno (2.13)
$$  
Correspondingly applying Eq.(2.11) to Eqs.(2.12) and (2.13) may lead to
Part (a) and
 Part (b) of this theorem.  \qquad $ \Box $ \\

\noindent {\Large {\bf 3. Two universal factorization equalities  on complex
quaternion matrices }} \\
  
In this section, we extend the  universal similarity factorization equality
in (2.4) to any $ m \times n $ matrix over $ \Bbb Q $, and give some of
 its consequences. 

\medskip
      
\noindent  {\bf Theorem 3.1.} \, {\em Let $ A = A_0 + A_1e_1 + A_2e_2 +
A_3 e_3 \in \Bbb Q^{ m \times n}$ be given. Then  $ A $ satisfies the
following universal factorization equality
$$
Q_{2m} \left[ \begin{array}{rr}  A  & 0   \\0 &  A \end{array} 
\right] Q_{2n}= \left[ \begin{array}{cc}  A_0 + A_1 i & -
( \, A_2 + A_3i \, ) \\ A_2 - A_3 i &
 A_0 - A_1 i  \end{array}
 \right] := \Psi(A) \in \Bbb C^{ 2m \times 2n},   \eqno (3.1) 
$$ 
where $Q_{2t}$ has the independent form 
$$ 
Q_{2t}= Q_{2t}^{-1} =Q^{\dagger}_{2t} = \frac{1}{2}
\left[ \begin{array}{cc} ( \,1 - ie_1 \,)I_t & (  \,e_2 + ie_3 \,)I_t \\
 ( \,-e_2 + i e_3 \,)I_t & ( \,1 + ie_1 \,)I_t   \end{array} \right], \
 \ \ \ \ t = m, \ n.  \eqno (3.2)
$$ 
In particular, when $ m = n,$ Eq.(3.1)becomes a universal similarity
factorization  equality over $ \Bbb Q$ }. 

\medskip

\noindent {\bf Proof.} \, It follows directly from multiplying out the
  three block
matrices  in the left-hand side of Eq.(3.1).    \qquad $ \Box$   

\medskip

The matrix $ \Psi(A)$ in Eq.(3.1) is called the complex representation of
 $ A $.  If setting
$$
 A_{11} = A_0 + A_1 i,  \ \ \ \ A_{12} =   -( \, A_2 + A_3i \, ), \ \ \ \
  A_{21} =  A_2 - A_3 i, \ \ \ \ A_{11} = A_0 - A_1 i
$$
in Eq.(3.1), then it can equivalently be expressed as  
$$
 Q_{2m}\left[ \begin{array}{rr}  A_{11}  & A_{12}  \\ A_{21} &  A_{22}
 \end{array}  \right]Q_{2n} = \left[ \begin{array}{rr}  A  & O  \\ 0 &
  A \end{array}  \right] \in \Bbb Q^{ 2m \times 2n}, \eqno(3.3)
$$ 
where $ A_{st} \in  \Bbb C^{ m \times n}$ are arbitrary, $ Q_{2t}$ is as
in Eq.(3.2), and  $ A $ has the following form
$$ 
A = \frac{1}{2}( \, A_{11} +  A_{22} \, ) +  \frac{1}{2}( \, A_{22} -
 A_{11} \, )ie_1 +  \frac{1}{2}( \, A_{21} -  A_{12} \, )e_2 +
 \frac{1}{2}( \, A_{12} +  A_{21} \, )ie_3  \in \Bbb Q^{ m \times n} .
  \eqno(3.4)
$$ 
which can also be stated that for any matrix  $M
\in \Bbb C^{ 2m \times 2n}$, there must be  a unique matrix $ A
 \in \Bbb Q^{ m \times n}$ such that
$$ 
\Psi(A) = M.  \eqno (3.5) 
$$ 

Various operation properties on complex representation of complex
quaternion matrices can easily be derived from Eq.(3.1). 

\medskip

\noindent {\bf Theorem 3.2.} \, {\em Let $ A, \,
B \in \Bbb Q^{ m \times n},  \,  C \in \Bbb Q^{ n \times p},$ and
$ \lambda \in \Bbb C $ be given.  Then

{\rm (a)} \ $ A=B \Longleftrightarrow  \Psi(A) = \Psi(B).$ 

{\rm (b)} \ $ \Psi(A + B)=  \Psi(A) + \Psi(B), \ \ \ \Psi(AC)= \Psi(A)
\Psi(C ),  \ \ \ \Psi( \lambda A)= \Psi(A \lambda ) = \lambda \Psi(A).$

{\rm (c)} \ $ \Psi(A^*) =  \left[ \begin{array}{cc}  O  & I_n   \\
-I_n &  O \end{array}
\right] \Psi^T(A) \left[ \begin{array}{cc}  O  & -I_m   \\ I_m &  O
\end{array}
\right].$ 

{\rm (d)} \  $ \Psi_1(A^{\dagger}) = \Psi^*(A),$ the conjugate transpose
of $\Psi(A).$

{\rm (e)} \ $ A = \frac{1}{4}E_{2m} \Psi(A)E_{2n}^{\dagger},$ where
$ E_{2t}= [\, ( \,1- ie_1 \,)I_t, \ ( \,e_2 + ie_3 \,)I_t \, ], \,
t = m , \, n .$

{\rm (f)} \ $ A$ is invertible if and only if $ \Psi(A) $ is invertible, 
in that case, $ \Psi(A^{-1}) = \Psi^{-1}(A)$ and $ A^{-1} =
\frac{1}{4}E_{2m} \Psi^{-1}(A)E_{2m}^{\dagger}.$

{\rm (g)}  \ $ \Psi(A)E_{2n}^{\dagger}E_{2n} = E_{2m}^{\dagger}E_{2m}
\Psi(A).$

{\rm (h)} \ $ A$ is Hermitian if and only if $ \Psi(A) $ is  Hermitian
over $\Bbb C$.

{\rm (i)} \ $ A$ is unitary if and only if $ \Psi(A) $ is  unitary over
 $\Bbb C$. } 

\medskip

Another universal factorization equality on complex quaternion matrices
is established as follows. 

\medskip

\noindent  {\bf Theorem 3.3.} \, {\em Let $ A = (a_{st})
\in \Bbb Q^{ m \times n}, $ and  denote
$ D_A = ( a_{st}I_2 ) \in \Bbb Q^{ 2m \times 2n}$.  Then $D_A$ satisfies
the following universal factorization equality
$$
P_{2m}D_AP_{2n}= \left[ \begin{array}{ccc} \psi(a_{11}) & \cdots &
 \psi(a_{1n}) \\ \vdots   &  & \vdots
 \\ \psi(a_{m1}) & \cdots &  \psi(a_{mn})   \end{array} \right] := \psi(A)
  \in \Bbb C^{ 2m \times 2n},   \eqno (3.6)
$$ 
where $P_{2t}$ has the form 
$$ 
P_{2t}= P^{-1}_{2t} = P^{\dagger}_{2t} = {\rm diag}( \, Q, \ \cdots, Q \, ),
 \ \ \ Q= Q^{-1} = Q^{\dagger} = \frac{1}{2} \left[ \begin{array}{cc} 1
  - ie_1 & e_2 + ie_3 \\ -e_2 + i e_3 & 1 + ie_1   \end{array} \right].
   \ \ \ \ \ t =m, \ n .
$$ 
In particular, when $ m = n,$ Eq.(3.6) becomes a universal similarity
factorization  equality over $ \Bbb Q.$ } 

\medskip

\noindent{\bf Proof.}\, Observe from Eq.(2.1) that $ Q(a_{st}I_2)Q =
\psi(a_{st})$. We immediately obtain
$$ 
P_{2m}D_AP_{2n} = [ \, Q(a_{st}I_2)Q \,]_{m \times n} =
[ \, \psi(a_{st}) \, ]_{m \times n} = \psi(A).
$$  
which is exactly Eq.(3.6).  \qquad $ \Box$ 

\medskip

The complex matrix $ \psi(A)$ in Eq.(3.6) is also called the complex
representation of $ A $. Clearly the two complex matrices
 $ \Psi(A)$ in Eq.(3.1)  and  $\psi(A)$ Eq.(3.6) are permutationally equivalent,
  that is,  there are two permutation matrices $ G$ and $H$ such that
  $ G \Psi(A)H = \psi(A)$. 

\medskip

For convenience of application, Eq.(3.6) can be simply stated that  for any
$ M \in \Bbb C^{ 2m \times 2n},$ there must exist an $ A \in
\Bbb Q^{m \times n}$ such that
$$ 
\psi(A) = M. \eqno (3.7)
$$      

\noindent {\bf Theorem 3.4.} \, { \em Let $ A, \,  B
\in \Bbb Q^{ m \times n},  \,  C \in \Bbb Q^{ n \times p}, \,
 \lambda \in \Bbb C $ be given.  Then

{\rm (a)} \ $ A=B\Longleftrightarrow  \psi(A) = \psi(B).$ 

{\rm (b)} \ $ \psi(A + B)=  \psi(A) + \psi(B), \ \ \ \psi(AC)= \psi(A)
\psi(C ),  \ \ \ \psi( \lambda A)= \psi(A \lambda ) = \lambda \psi(A).$

{\rm (c)} \ $\psi(I_m) = I_{2m}.$ 

{\rm (d)} \ $ \psi(A^{\dagger}) = \psi^*(A),$ the conjugate transpose of
$\psi(A).$

{\rm (e)} \ $ A$ is invertible if and only if $ \psi(A) $ is invertible,
in that case, $ \psi(A^{-1}) = \psi^{-1}(A).$

{\rm (h)} \ $ A$ is Hermitian if and only if $ \psi(A) $ is  Hermitian.  

{\rm (i)} \ $ A$ is unitary if and only if $ \psi(A) $ is  unitary. } 

\medskip

Just as for  complex  quaternions, we can define the Moore-Penrose
inverse of any $ m \times n $ complex quaternion matrix as follows.

\medskip
 
\noindent {\bf Definition.} \, Let $ A \in \Bbb Q^{ m \times n} $
be given. If the following four equations
$$ 
AXA  = A,  \qquad  XAX  = X, \qquad  ( AX )^{\dagger} = AX,
 \qquad   ( XA )^{\dagger} = XA \eqno (3.8)
$$ 
have a common solution for $ X  $, then  this solution is called the
Moore-Penrose inverse of $ A $, and denoted by $ X = A^{+}$. 

\medskip

The existence and the uniqueness of the Moore-Penrose inverse of a
 complex quaternion matrix $ A $  can be determined by its complex
 representation $ \Psi( A)$. In fact, according to Theorem 3.2(a),
 (c) and (d), the four equations in Eq.(3.8) are equivalent to the
  following four  complex matrix equations
$$  
 \Psi(A)\Psi(X)\Psi(A) = \Psi(A), \qquad 
 \Psi(X)\Psi(A)\Psi(X) = \Psi(X),
$$
$$
[ \, \Psi(A)\Psi(X) \, ]^* = \Psi(A) \Psi(X) , \qquad 
[\, \Psi( X)\Psi(A) \,]^* = \Psi(X) \Psi(A).
 $$ 
According to the complex matrix theory, the following four equations  
$$  
 \Psi( A)Y \Psi( A) = \Psi( A), \qquad    Y \Psi(A)Y = Y,
  \qquad  [ \, \Psi(A)Y \,]^* = \Psi(A)Y , \qquad 
   [ \, Y \Psi(A) \,]^* = Y  \Psi(A)
$$ 
have  a unique common solution $  Y = \Psi^{+}(A)$, the Moore-Penrose
of $ \Psi(A)$. Then by Eq.(3.5), there is a unique matrix $ X $ over
$ \Bbb Q $ such that $   \Psi( X ) = Y =  \Psi^{+}(A)$,  in which case,
this $ X $ can be expressed  as
$$
 X = \frac{1}{4}E_{2n}YE_{2m}^{\dagger} = \frac{1}{4} E_{2n}
 \Psi^{+}(A)E_{2m}^{\dagger}.
$$ 
Correspondingly  this $ X $ is the unique solution to  Eq.(3.8). Hence  we
have  the following.

\medskip

\noindent {\bf Theorem 3.5.} \, { \em Let $ A \in \Bbb Q^{ m \times n}$
be given. Then its Moore-Penrose inverse  $ A^{+} $ of $ A $ exists
uniquely, and  satisfies the following two equalities
$$ 
 \Psi(A^{+}) = \Psi^{+}(A), \qquad   A^{+} = \frac{1}{4} E_{2n}
 \Psi^{+}(A)E_{2m}^{\dagger},
$$ 
where $E_{2t}= [ \, (\, 1 - ie_1 \,)I_t, \ ( \,e_2 +  ie_3 \,)I_t \, ],
 \, t = m, \,  n .$} 

\medskip

Through the complex representation in Eq.(3.1), we can define the rank of any
 $ A \in \Bbb Q^{m \times n}$ as follows $ {\rm rank}( A )
  = \frac{1}{2}{\rm rank}[ \, \Psi(A) \,].$  Obviously the rank of a
  complex quaternion matrix $ A $ is a fraction if the rank of $ \Psi(A) $
   is an odd number. In particular, $ A \in \Bbb Q^{m \times m }$ is
    invertible if and only if  $ {\rm rank}( A )= m.$ 

\medskip
  
On the basis of the above results  we now are able to investigate various
kinds of  problems related to complex  quaternion matrices.  In the next
 three sections, we  consider three basic problems---right eigenvalues and
 eigenvectors, similarity, as well as determinants of square complex
 quaternion matrices. \\

\noindent {\Large {\bf 4. Right eigenvalues and eigenvectors of complex quaternion
 matrices }}\\

As usual, {\em right eigenvalue equation} for a complex quaternion matrix
 $ A \in  \Bbb Q^{ n\times n} $ is defined by
$$ 
AX = X\lambda,  \ \ \  \ X  \in  \Bbb Q^{n \times 1}, \ \ \ \  \lambda
  \in  \Bbb Q.  \eqno (4.1)
$$
If a  $ \lambda  \in  \Bbb Q $ and a nonzero $  X  \in
\Bbb Q^{n \times 1}$ satisfy Eq.(4.1), then   $ \lambda $ is called a
{\em right eigenvalue} of $ A $ and $ X $ is called  an {\em eigenvector}
associated with $ \lambda $. In particular, if a $ \lambda \in  \Bbb Q$
and  an $ X  \in  \Bbb Q^{n \times 1}$ with rank$(X) = 1 $ satisfy Eq.(4.1),
then $ \lambda  $ is called a {\em regular right eigenvalue} of $ A $ and $ X $ is called  an {\em regular eigenvector} associated with $ \lambda $. In this section, we shall prove that any square complex quaternion matrix has at least  one complex right eigenvalue, and  also has at least one regular right eigenvalue. To do so, we need some preparation. \\

\noindent {\bf Definition.} \,  For any $ X = X_0 + X_1e_1 + X_2e_2 + X_3 e_3
 \in \Bbb Q^{ n \times 1} $, we call the following complex matrix
$$  
\left[ \begin{array}{r}  X_0 + X_1i    \\  X_2 - X_3i \end{array}  \right]
:= \overrightarrow{X }
 \in \Bbb C^{ 2n \times 1},  \eqno (4.2) 
$$ 
uniquely determined by $ X $, the {\em complex adjoint vector} of $ X $.

\medskip

\noindent {\bf Lemma 4.1.} \, {\em  Let $A \in \Bbb Q^{ m \times n }, \ X
  \in \Bbb Q^{ n \times 1 },$ and $ \lambda \in \Bbb C$ be given. Then
  $ \overrightarrow{AX } = \Psi(A) \overrightarrow{X } $ and  $
  \overrightarrow{X \lambda } = \overrightarrow{X}  \lambda. $ } 

\medskip

{\bf Proof.} \, It is easy to see that for any column matrix $ Y
\in \Bbb Q^{ n \times 1 } $, the corresponding $
 \Psi(Y) $  and  $\overrightarrow{Y} $ satisfy the relation
 $ \overrightarrow{Y } =  \Psi(Y)[ \, 1, \ 0 \, ]^T. $ Now applying it to
  $ AX $ and $ X \lambda$ we find
$$
\overrightarrow{AX} =  \Psi(AX )[ \, 1,  \ 0 \, ]^T =  \Psi(A)
\Psi(X)[ \, 1, \ 0 \, ]^T = \Psi(A) \overrightarrow{X },
$$
$$
 \overrightarrow{X \lambda } =  \Psi(X \lambda )[ \, 1, \  0 \, ]^T =
  \Psi(X) \lambda [ \, 1,  \ 0 \, ]^T
 = \overrightarrow{X } \lambda. \qquad \Box 
$$ 

Based on the above notation, we can deduce the following several results.

\medskip

\noindent {\bf Theorem 4.2.} \, {\em  Let $A \in \Bbb Q^{ n \times n }$
be given. Then all the eigenvalues of $ \Psi(A)$ are  right eigenvalues of
 $A ,$ and all the  eigenvectors of $  \Psi(A ) $ can be used for
 constructing eigenvectors of $ A. $ } 

\medskip

\noindent {\bf Proof.} \, Assume that $ \lambda \in \Bbb C$ and $ O
\neq Y \in \Bbb C^{ 2n \times 1 }$ satisfy
$$
\Psi(A)Y = Y \lambda. 
$$ 
Then we set $ X = E_{2n} Y \in \Bbb Q^{ n \times 1 },$ where $
 E_{2n} = [ \,  ( \, 1- ie_1 \, )I_n, \   ( \,  e_2 + ie_3 \, )I_n \, ] $,
 and it satisfies $  E_{2n} E_{2n}^{\dagger} = 4I_n.$ Combining the above
 results with Theorem 3.2(e) and (g), we obtain
$$ 
AX  =    \frac{1}{4}E_{2n} \Psi(A)E_{2n}^{\dagger}E_{2n}Y  =
  \frac{1}{4}E_{2n} E_{2n}^{\dagger}E_{2n} \Psi(A)Y  =  E_{2n} Y
  \lambda = X \lambda,
$$ 
which shows that $ \lambda$ is a right eigenvalue of $ A $ and $ X =
  E_{2n}Y $ is an eigenvector of $ A$ associated with this $ \lambda$.
  So the conclusion of the theorem is true.  \qquad $ \Box $ 

\medskip
 
Conversely, we have the following result. 

\medskip

\noindent {\bf Theorem 4.3.} \, {\em  Let $A \in \Bbb Q^{ n \times n }$ be
 given. Then all the complex right eigenvalues of $A$ are eigenvalues of $
 \Psi(A ),$ too.} 

\medskip

\noindent  {\bf Proof.} \, Assume that  $ \lambda \in \Bbb C$ and $ O
 \neq X  \in \Bbb Q^{ n \times 1 }$ satisfy $ AX = X \lambda. $ Then
 applying Lemma 4.1 to the both sides of this equality, we find $ \Psi(A)
  \overrightarrow{X} =  \overrightarrow{X} \lambda, $ which shows that
   $ \lambda $ and $  \overrightarrow{X}$ are a pair of eigenvalue and
   eigenvector of  $  \Psi(A) $.  \qquad $ \Box $ 

\medskip

Next are two results on the regular right eigenvalues and eigenvectors
of complex quaternion matrices. 

\medskip

\noindent {\bf Theorem 4.4.} \, {\em  Let $A \in \Bbb Q^{ n \times n }$
be given. Then $ A $ has at least  one regular right eigenvalue.} 

\medskip

\noindent {\bf Proof.}  \, We first assume that the complex  representation
 $ \Psi(A) $ has at least two linearly independent  eigenvectors  $Y_1 $
 and $ Y_2$, i.e.,
$$    
 \Psi(A)Y_1  = Y_1\lambda_1, \ \ \ \  \Psi(A)Y_2  = Y_2 \lambda_2,
 \eqno (4.3)
$$
where ${\rm rank}[ \, Y_1, \ Y_2 \, ] = 2,$  $\lambda_1 $ and $ \lambda_2$
are two complex numbers with  $ \lambda_1 = \lambda_2 $, or
$ \lambda_1 \neq \lambda_2 $.
Then 
$$
\Psi(A)[ \, Y_1, \ Y_2 \, ] = [ \, Y_1,  \ Y_2 \, ] \left[ \begin{array}{rr}
 \lambda_1   & 0   \\ 0 &  \lambda_2 \end{array}
\right].  \eqno (4.4) 
$$ 
According to Eq.(3.5),  there must exist an $X \in \Bbb Q^{ n \times 1 }$
such that
$$ \Psi(X) = [ \, Y_1, \ Y_2 \, ], \ \ \ \ {\rm rank}( X ) = \frac{1}{2}
{\rm rank}\Psi(X) = 1, \eqno (4.5)
$$ 
meanwhile there must exist a $ \lambda \in \Bbb Q$ such that  
$$
 \Psi( \lambda )=  \left[ \begin{array}{rr} \lambda_1   & 0   \\ 0 &
 \lambda_2 \end{array}
\right].  \eqno (4.6) 
$$ 
Putting Eqs.(4.5) and (4.6) into Eq.(4.4) and then applying Theorem 3.2(a)
to it, we obtain
$$
\Psi(A) \Psi(X) = \Psi(X)\Psi( \lambda ) \Longrightarrow  AX = X \lambda,
 \ \ and \ \ {\rm rank}( X) = 1,
$$
which shows that  $ \lambda $ and $ X $ are a pair of regular right
eigenvalue and eigenvector of $ A $. 

\medskip

Next assume that the complex  matrix $ \Psi(A) $  has only one  eigenvalue
$ \lambda_1$ and only one linearly independent eigenvector $ Y_1$
associated with  $ \lambda_1 $. Then according to  complex matrix theory,
there exists another vector $ Y_2$  over $ \Bbb C $ such that
$$ 
\Psi(A)[\, Y_1, \ Y_2 \, ] = [ \, Y_1, \ Y_2 \, ] \left[ \begin{array}{cc}
 \lambda_1   & 1   \\ 0  &  \lambda_1 \end{array}   \right], \eqno (4.7)
$$ 
where rank$[ \, Y_1, \ Y_2 \, ] = 2. $ Then it follows by Eq.(3.5) that  there
must exist a  $ \lambda \in \Bbb Q$ and an $X \in \Bbb Q^{ n \times 1 }$
such that
$$
  \Psi( \lambda )=  \left[ \begin{array}{cc} \lambda_1   & 1   \\ 0 &
   \lambda_1 \end{array} \right], \ \ \Psi(X) = [ \, Y_1, \ Y_2 \, ],
   \ \ \ \ {\rm rank}( X ) = \frac{1}{2} {\rm rank} \Psi(X) = 1.
    \eqno (4.8)
$$ 
In that case, combining Eqs.(4.7) and (4.8) and Theorem 3.2(a), we get 
$$
\Psi(A) \Psi(X) = \Psi(X)\Psi( \lambda ) \ \Longrightarrow \ AX = X \lambda,
 \ \ r( X ) = 1,
$$
which shows that  $ \lambda $ and $ X $ are a pair of regular right
 eigenvalue and eigenvector of  $ A $. \qquad $ \Box $ 

 \medskip
 
\noindent {\bf Theorem 4.5.} \, {\em  Let $A \in \Bbb Q^{n \times n }$ be
given. Then from any regular right eigenvalue of $A$ we can derive
eigenvalues of $ \Psi(A). $ } \\

\noindent {\bf Proof.}  \ Assume that  $ \lambda =   \lambda_0 + \lambda_
1e_1 + \lambda_2e_2 + \lambda_3e_3  \in \Bbb Q$ and $ X
 \in \Bbb Q^{ n \times 1 }$ with $ {\rm rank}( X ) = 1 $  satisfy
$$
AX = X \lambda. \eqno (4.9)
$$ 
If $  \lambda \in \Bbb C$ in Eq.(4.9), then from Theorem 3.3(a) and (b)
 we get
$$
\Psi(A) \Psi(X) = \Psi(X)\Psi( \lambda ) = \Psi(X)
\left[ \begin{array}{rr} \lambda & 0   \\ 0 &  \lambda \end{array} \right],
$$ 
which clearly  shows that $ \lambda $ is an  eigenvalue  of $ \Psi(A)$, and
  $ \Psi(X)$ are two eigenvectors of $ \Psi(A)$. If $  \lambda
  \notin \Bbb C$ in Eq.(4.9), then by Theorem 2.6 we know that  this
   $ \lambda $ can be expressed as
$$  
\lambda = p[\,  \lambda_0 + \tau(\lambda ) e_1 \, ]p^{-1}, \eqno (4.10)  
$$ 
when $ \tau^2(\lambda )  = \lambda_1^2 + \lambda_2^2 + \lambda_3^2 \neq 0 $,
 or
$$  
    \lambda = q \left( \,  \lambda_0 - \frac{1}{2} e_2 + \frac{1}{2} ie_3 \,
     \right)q^{-1}, \eqno (4.11)
$$ 
when $ \lambda_1^2 + \lambda_2^2 + \lambda_3^2 = 0 $. Under the condition in
(4.10), the equality in Eq.(4.9) can be equivalently expressed as
$$  
A(Xp) =(Xp)[ \,  \lambda_0 + \tau(\lambda ) e_1 \, ]. \eqno(4.12) 
$$
Applying Theorem 3.2(b) to the both sides of Eq.(4.12) yields 
$$
\Psi(A) \Psi(Xp) =  \Psi_1(Xp)\Psi[ \,  \lambda_0 + \tau(\lambda )e_1 \, ]
= \Psi(Xp) \left[ \begin{array}{cc} \lambda_0 + \tau(\lambda )i  & 0
  \\ 0 &  \lambda_0 - \tau(\lambda )i \end{array} \right],
$$ 
which  shows that $ \lambda_0 \pm  \tau(\lambda )i $ are two eigenvalues
 of $ \Psi(A)$, and  $ \Psi(Xp)$ are two eigenvectors of $ \Psi(A)$.  Under
 the condition in Eq.(4.11), the equality in (4.9) can  equivalently be
 expressed as
$$  
A(Xq) =(Xq)\left( \,  \lambda_0 - \frac{1}{2} e_2 + \frac{1}{2} ie_3 \,
\right).   \eqno(4.13)
$$ 
Then applying Theorem 3.2(b) to the both sides of Eq.(4.13) yields 
$$
\Psi(A) \Psi(Xq) = \Psi(Xq)\Psi\left(  \, \lambda_0 - \frac{1}{2} e_2 +
\frac{1}{2} ie_3 \, \right) = \Psi(Xq) \left[ \begin{array}{cc} \lambda_0
 & 1   \\ 0 &  \lambda_0  \end{array} \right],
$$ 
which  shows that $ \lambda_0 $ is an eigenvalue  of $ \Psi(A)$, and the
first column of   $ \Psi(Xp)$ is  an eigenvector of $ \Psi(A)$.
\qquad $ \Box $ \\

\noindent {\Large {\bf 5. Similarity  of complex quaternion matrices }} \\ 

Two square matrices $ A $ and $ B $ of the same size over $ \Bbb Q $
are said to be {\em similar} if there exists an invertible matrix $ X $
over $ \Bbb Q $  such that $ X^{-1}A X = B $. Based on the results in
the preceding two subsections, we can easily find a simple result
charactering the similarity of two complex quaternion matrices. 

\medskip

\noindent {\bf Theorem 5.1.} \, {\em  Let $A, \, B \in
\Bbb Q^{n \times n}$ be given. Then the following three statements
  are equivalent:

{\rm (a)} \ $ A $ and $ B $ are similar  over $\Bbb Q.$  

{\rm (b)} \ $\Psi(A) $ and  $ \Psi(B)$ are  similar  over $\Bbb C$.

{\rm (c)} \ $\psi(A) $ and  $ \psi(B)$ are  similar  over $\Bbb C$. } 

\medskip 
 
\noindent {\bf Proof.} \,  Assume first that $ A \sim B $ over $\Bbb Q$.
Then there is an invertible matrix $ X $ such that $X^{-1}A X = B$. Now
applying Theorem 3.2(b) and (f) to its both sides yields  $\Psi^{-1}( X)
\Psi( A ) \Psi( X ) = \Psi(B)$, which shows that $ \Psi(A) \sim  \Psi(B)
 $ over $\Bbb C$. Conversely, assume that $ \Psi(A) \sim  \Psi(B)$ over $\Bbb C$. Then there is an invertible matrix $ Y \in  \Bbb C^{ 2n \times 2n}$ such that $Y^{-1}\Psi(A)Y =\Psi(B)$. For this $ Y $, by (3.5) there must be an invertible matrix $ X \in  \Bbb Q^{ n \times n }$ such that $ \Psi(X) = Y $. Thus, $\Psi( X^{-1})\Psi(A) \Psi(X) =\Psi(B)$, and consequently $ X^{-1}AX = B $, which shows that  $ A \sim B $. The equivalence of (a) and (c) can also be shown in the same manner.
 \qquad  $ \Box $ 

\medskip

Based on this result, we can extend various results on similarity of complex
matrices to complex quaternion matrices. 

\medskip 

\noindent {\bf Theorem 5.2.} \, {\em  Let $A \in \Bbb Q^{ n \times n }$ be
given. Then  $ A $ is similar to a diagonal matrix  over $\Bbb Q$ if and
only if the sizes of the Jordan blocks in the Jordan
canonical form of  the complex representation $\psi(A)$ of $ A $  are not
greater than 2. } 

\medskip 

\noindent {\bf Proof.} \, Assume that $ A $ is diagonalizable over
$\Bbb Q$, i.e.,  there is an invertible matrix $ P $ over $\Bbb Q$ such
that
$$  
P^{-1}AP = {\rm diag}( \, \lambda_1, \ \lambda_2, \  \cdots, \
 \lambda_n \,),  \eqno(5.1)
$$ 
where $\lambda_1$---$\lambda_n \in \Bbb Q$. Then by Theorem 5.1, we obtain 
 $$  
\psi^{-1}(P)\psi(A)\psi(P) = {\rm diag}( \,\psi( \lambda_1), \
 \psi( \lambda_2), \  \cdots, \ \psi( \lambda_n)  \, ),  \eqno(5.2)
$$   
where $\psi( \lambda_t)$ is a $ 2 \times 2 $ complex matrix. This equality
implies that the sizes of Jordan block  in the Jordan canonical form of
 $\psi(A)$  are not greater than 2. Conversely assume that
there is an invertible  matrix $ G$ over $ \Bbb C$ such that 
$$  
G^{-1}\psi(A)G = {\rm diag}( \, J(\mu_1), \ \cdots, \ J(\mu_k), \
 \mu_{k+1}, \  \cdots, \ \,  \mu_{r} \,),  \eqno(5.3)
$$   
where $ J(\mu_t)(t = 1$---$k$) is a $2 \times 2$ complex Jordan block,
 $ \mu_{k+1}$---$\mu_{r}$ are
complex eigenvalues of $\psi(A)$ with $2k + ( r - k ) = 2n$. In that case,
 by Eq.(3.7) we know that there is an invertible matrix $ P $ over
 $ \Bbb Q$ such that $\phi(P) = G$, there is a $ \lambda_t$ over
$\Bbb Q$ such that $ \psi( \lambda_t) =J(\mu_t)(t = 1$---$k$), and there
are $ \lambda_{k+1}, \cdots, \
\lambda_{r} \in \Bbb Q$ such that 
$$
 \psi(\lambda_{k+1}) = \left[ \begin{array}{cc} \mu_{k+1}  & 0   \\ 0 &
 \mu_{k+2}\end{array} \right], \ \  \psi(\lambda_{k+2}) =
  \left[ \begin{array}{cc} \mu_{k+3}  & 0   \\ 0 & \mu_{k+4}\end{array}
   \right],
\ \cdots, \ \ \psi(\lambda_{n}) = \left[ \begin{array}{cc} \mu_{r-1}  & 0
 \\ 0 & \mu_{r}\end{array} \right], \eqno (5.4)
$$ 
In that case, let $ D = {\rm diag}( \, \lambda_1, \   \lambda_2, \  \cdots,
\  \lambda_n \,)$. Then Eq.(5.3) becomes
$$  
\psi^{-1}(P)\psi(A)\psi(P) = \psi(D) \ \Longrightarrow \  P^{-1}AP = D. 
$$
Hence $ A $ is diagonalizable over $ \Bbb Q$. \qquad $\Box$. 

\medskip 

The result in the above theorem alternatively implies that if the size of a
 Jordan block in the Jordan
canonical form of the complex representation $\psi_(A)$ of $ A $ is greater
 than 2, then $ A $ can not be diagonalizable over $\Bbb Q$. 

\medskip 

\noindent {\bf Theorem 5.3.} \, {\em  Let $A \in \Bbb Q^{ n \times n }$ be
given. Then  $ A $ is similar to a complex matrix over $\Bbb Q$ if and
only if  the complex representation $\Psi(A)$ of $ A $  is similar to a
block diagonal matrix ${\rm diag}(\, J, \ J \, )$ over $\Bbb C,$ where
$ J \in \Bbb C^{ n \times n }$. } 

\medskip 

\noindent {\bf Proof.} \, If $ A \sim J \in \Bbb C^{ n \times n }$, then
it follows  by Theorem 5.1 and Eq.(3.1) that $ \Psi(A) \sim \Psi(J)
= {\rm diag}( \, J, \ J \,)$. Conversely if $ \Psi(A)
\sim {\rm diag}( \, J, \ J \,)$ with $ J \in \Bbb C^{ n \times n }$,
then it follows by (3.1) that $ \Psi(J) = {\rm diag}( \, J, \ J \,)$.
Consequently $ A \sim J$ follows by Theorem 5.1.  \qquad $\Box$ \\

\noindent {\Large {\bf 6. Determinants of  complex quaternion matrices }} \\ 

As one the most fundamental problems in the theory of quaternion matrices,
the  determinants of complex quaternion matrices, including the
determinants of real quaternion matrices, have been considered  by lots
of authors from different aspects, see, e.g., \cite{3}, \cite{5}, \cite{6}
and \cite{8}. Among all of them, a direct and simple method for defining
determinants of quaternion matrices is through their representations in the
central fields of the corresponding quaternion algebras. In \cite{8},
Zhang presents this kind of definition for determinants of real
quaternion matrices, and demonstrates the consistency of his definition
 with some other classic definitions on determinants based on the products
 and  sums of entries in matrices.  

\medskip

Following the introduction of the universal  similarity equalities for
generalized quaternion matrices, we shall easily find that it is a
most reasonable method to define determinants of quaternion matrices
through their representations in the central field of  corresponding
quaternion algebra. In fact, from Eq.(3.1), we know that for any $ A
\in \Bbb Q^{ n \times n }$
$$
Q_{2n} {\rm diag}( \,  A, \, A  \, ) Q^{-1}_{2n}= \Psi(A)
\in \Bbb C^{ 2n \times 2n}.   \eqno (6.1)
$$
Therefore if we hope that determinants of quaternion matrices
over $\Bbb Q $ satisfy the following two basic properties
$$ 
{\rm det}( MN )={\rm det} \, M{\rm det} N  \ \ \ \ and   \ \ \ \
{\rm det} I_n = 1
$$   
for any $ M, \ N \in \Bbb Q^{ n \times n}$,  then as a natural
consequence of these two properties, the determinant of the diagonal
matrix $ {\rm diag}( \, A, \, A \, )$  satisfies the following equality
$$
 {\rm det}[ \,{\rm diag}( \, A, \, A \, ) \, ] ={\rm det}[ \, Q_{2n}
  {\rm diag}(  \, A, \, A \, ) Q^{-1}_{2n} \, ] = {\rm det} \Psi(A) .
  \eqno (6.2)
$$ 
From this equality we see that there seem only two kinds of natural
choices for the definition of determinant of $ A \in
\Bbb Q^{ n \times n}$, namely, define
$$
{\rm det}A := | \Psi(A) |,   \eqno (6.3) 
$$ 
or
$$  
{\rm det}A := | \Psi(A) |^{\frac{1}{2}},   \eqno (6.4) 
$$  
where $| \Psi(A)| $  is the ordinary determinant of the complex matrix
$ \Psi(A) $. For simplicity, we prefer Eq.(6.3) as the definition of
determinant of matrix over $\Bbb Q $ and call it the {\em central
determinant} of matrix $A$, and  denote it by $ | A|_c $.  

\medskip

Some basic operation properties on the central determinants of complex
quaternion matrices are listed below without proofs. 

\medskip 

\noindent {\bf Theorem 6.1.}  \, {\em  Let $A, \, B  \in
\Bbb Q^{ n \times n }, \, \lambda \in \Bbb C,$ and $ \mu \in
 \Bbb Q$ be given.

{\rm (a)} \ If $ A  \in \Bbb C^{ n \times n},$ then $
 |A|_c = |A|^2. $

{\rm (b)} \ $ A $  is invertible  $ \Longleftrightarrow$ 
 $ |A|_c \neq 0 $.

{\rm (c)} \ $ | AB|_c = |A|_c \,|B|_c.$

{\rm (d)} \  $ |\lambda A|_c = |A \lambda |_c = \lambda^{2n} |A|_c. $ 
 
{\rm (e)}  \ $ |\mu A|_c = |A \mu |_c = n^{2n}( \mu ) |A|_c,$
 where  $ n( \mu ) $ is the weak norm of $ \mu $.

{\rm (f)} \  $  |A^{-1}|_c = |A|^{-1}_c. $ 
 
{\rm (g)} \  $ |A^{\dagger}|_c = \overline{ |A|}_c.$ 
 
 {\rm (h)} \ If $ A = \left[ \begin{array}{rrc} a_{11} & \cdots & *
  \\  & \ddots & \vdots \\  & &   a_{nn} \end{array} \right], $ then
   $  |A|_c = n(a_{11}) n(a_{22}) \cdots  n(a_{nn}). $ \\
  
{\rm (i)} \ If $ A = \left[ \begin{array}{cc} A_1   & *   \\ 0 & A_2
 \end{array} \right]$, where $  A_1 $ and $ A_2$
are square, then $|A|_c =|A_1|_c  |A_2|_c $.  
 
 {\rm (j)} \ If $ A \sim B ,$ then $  |A|_c =  |B|_c. $ }

\medskip 
 
By means of the determinants of complex quaternion matrices defined above,
we can easily define the central  characteristic
polynomial of any $A \in \Bbb Q^{ n \times n }$ as follows
$$ 
 p_A( \lambda ) = | \lambda I_{2n} - \Psi(A)|, \eqno (6.5) 
$$
which is a complex polynomial of degree $ 2n$. From it we easily get the
following. 

\medskip

\noindent {\bf Theorem 6.2}(Cayley-Hamilton theorem).
\ {\em  Let $A \in \Bbb Q^{ n \times n }$ be given.
 Then $ p_A( A ) = 0 $. } 

\medskip

\noindent {\bf Proof.}  \, Since $ p_A( \lambda ) =
| \lambda I_{2n} - \Psi(A)|$, so  $ p_A[\Psi(A) ] = 0.$ Putting Eq.(6.1) in
it and simplifying the equality yields
$$ 
  p_A[ \, {\rm diag}( A, \, A ) \, ] = {\rm diag} ( \,  p_A(A), \
   p_A(A) \, )  = 0,
$$ 
which implies that $ p_A( A ) = 0 $.   \quad $ \Box $  \\

\noindent{\Large  {\bf 7. Conclusions}} \\ 

In the article, we have established a fundamental universal similarity
factorization equality over the complex quaternion algebra $ \Bbb Q$. This
equality clearly reveals the intrinsic relationship between the complex
quaternion algebra $ \Bbb Q$ and the $2 \times 2$ total matrix algebra,
and could serve as a powerful tool for investigating various problems related
to complex quaternions and their applications. In addition to the results in
Sections 3---6, we can also apply Eqs.(2.4), (3.1) and (3.6) to investigate
some other basic topics related to complex quaternion matrices, such as,
singular value decompositions, norms, numerical ranges of complex
quaternion matrices and so on. Finally we point out that the equality
(2.4) can also extended to the complex Clifford algebra $ \Bbb C_n$,
and a set of perfect matrix representation theory on the complex Clifford
algebra $ \Bbb C_n$ can explicitly established. We shall examine  this 
topic in  another paper. \\

\end{document}